\documentclass{amsart}
\usepackage{amsthm}
\usepackage{amstext}
\usepackage{amssymb}
\usepackage{esint}

\makeatletter
\numberwithin{equation}{section}
\numberwithin{figure}{section}
\theoremstyle{plain}
\newtheorem{thm}{Theorem}

\theoremstyle{remark}
\newtheorem*{rmk}{Remark}
\makeatother

\allowdisplaybreaks
\usepackage[sort&compress,square,numbers]{natbib}

\begin{document}

\title[degenerate poly-Bernoulli numbers and polynomials]{Some applications of degenerate poly-Bernoulli numbers and polynomials }
\author{Dae San Kim}
\address{Department of Mathematics, Sogang University, Seoul 121-742, Republic
of Korea}
\email{dskim@sogang.ac.kr}

\author{Taekyun Kim}
\address{Department of Mathematics, Kwangwoon University, Seoul 139-701, Republic
of Korea}
\email{tkkim@kw.ac.kr}

\begin{abstract}
In this paper, we consider degenerate poly-Bernoulli numbers and polynomials
associated with polylogarithmic function and $p$-adic invariant integral
on $\mathbb{Z}_{p}$. By using umbral calculus, we derive  some identities of those numbers
and polynomials.
\end{abstract}

\maketitle
\global\long\def\relphantom#1{\mathrel{\phantom{{#1}}}}

\global\long\def\Zp{\mathbb{Z}_{p}}

\global\long\def\li{\mathrm{Li}}

\global\long\def\acl#1#2{\left\langle \left.#1\right|#2\right\rangle }

\global\long\def\acr#1#2{\left\langle #1\left|#2\right.\right\rangle }

\section{Introduction}

Let $p$ be a fixed prime number. Throughout this paper, $\Zp$, $\mathbb{Q}_{p}$
and $\mathbb{C}_{p}$ will denote the ring of $p$-adic integers,
the field of $p$-adic rational numbers and the completion of the algebraic
closure of $\mathbb{Q}_{p}$. The $p$-adic norm is normalized
as $\left|p\right|_{p}=\frac{1}{p}$. For $k\in\mathbb{Z}$, the polylogarithmic
function $\li_{k}\left(x\right)$ is defined by $\li_{k}\left(x\right)=\sum_{n=1}^{\infty}\frac{x^{n}}{n^{k}}$.
For $k=1$, we have $\li_{1}\left(x\right)=-\log\left(1-x\right)$.

In \cite{key-3}, L. Carlitz considered the degenerate Bernoulli polynomials
which are given by the generating function
\begin{equation}
\frac{t}{\left(1+\lambda t\right)^{\frac{1}{\lambda}}-1}\left(1+\lambda t\right)^{\frac{x}{\lambda}}=\sum_{n=0}^{\infty}\beta_{n,\lambda}\left(x\right)\frac{t^{n}}{n!}.\label{eq:1}
\end{equation}

Note that $\lim_{\lambda\rightarrow0}\beta_{n,\lambda}\left(x\right)=B_{n}\left(x\right)$,
where $B_{n}\left(x\right)$ are the ordinary Bernoulli polynomials.
When $x=0$, $\beta_{n,\lambda}=\beta_{n,\lambda}\left(0\right)$
are called the degenerate Bernoulli numbers.

It is known that the poly-Bernoulli polynomials are defined by the
generating function
\begin{equation}
\frac{\li_{k}\left(1-e^{-t}\right)}{e^{t}-1}e^{xt}=\sum_{n=0}^{\infty}B_{n}^{\left(k\right)}\left(x\right)\frac{t^{n}}{n!},\quad\left(\text{see \cite{key-9}}\right).\label{eq:2}
\end{equation}
When $x=0$, $B_{n}^{\left(k\right)}=B_{n}^{\left(k\right)}\left(0\right)$
are called the poly-Bernoulli numbers.

Let $UD\left(\Zp\right)$ be the space of uniformly differentiable
functions on $\Zp$. For $f\in UD\left(\Zp\right)$, the $p$-adic
invariant integral on $\Zp$ is defined by
\begin{align}
\int_{\Zp}f\left(x\right)d\mu_{0}\left(x\right) & =\lim_{N\rightarrow\infty}\sum_{x=0}^{p^{N}-1}f\left(x\right)\mu_{0}\left(x+p^{N}\Zp\right)\label{eq:3}\\
 & =\lim_{N\rightarrow\infty}\frac{1}{p^{N}}\sum_{x=0}^{p^{N}-1}f\left(x\right),\quad\left(\text{see \cite{key-13}}\right).\nonumber
\end{align}

From (\ref{eq:3}), we have
\begin{equation}
\int_{\Zp}f\left(x+1\right)d\mu_{0}\left(x\right)-\int_{\Zp}f\left(x\right)d\mu_{0}\left(x\right)=f^{\prime}\left(0\right),\label{eq:4}
\end{equation}
where $f^{\prime}\left(0\right)=\left.\frac{df\left(x\right)}{dx}\right|_{x=0}$
(see \cite{key-1,key-2,key-3,key-4,key-5,key-6,key-7,key-8,key-9,key-10,key-11,key-12,key-13,key-14,key-15,key-16,key-17}).

By (\ref{eq:4}), we get
\begin{align}
\int_{\Zp}\left(1+\lambda t\right)^{\left(x+y\right)/\lambda}d\mu_{0}\left(y\right) & =\frac{\log\left(1+\lambda t\right)^{\frac{1}{\lambda}}}{\left(1+\lambda t\right)^{\frac{1}{\lambda}}-1}\left(1+\lambda t\right)^{\frac{x}{\lambda}}\label{eq:5}\\
 & =\frac{\log\left(1+\lambda t\right)}{\lambda t}\frac{t}{\left(1+\lambda t\right)^{\frac{1}{\lambda}}-1}\left(1+\lambda t\right)^{\frac{x}{\lambda}}\nonumber \\
 & =\sum_{n=0}^{\infty}\left(\sum_{l=0}^{n}\binom{n}{l}\lambda^{n-l}D_{n-l}\beta_{l,\lambda}\left(x\right)\right)\frac{t^{n}}{n!},\nonumber
\end{align}
where $D_{n}$ are the Daehee numbers of the first kind
given by the generating function
\begin{equation}
\frac{\log\left(1+t\right)}{t}=\sum_{n=0}^{\infty}D_{n}\frac{t^{n}}{n!},\quad\left(\text{see \cite{key-10}}\right).\label{eq:6}
\end{equation}

Let $\mathcal{F}=\left\{ \left.f\left(t\right)=\sum_{k=0}^{\infty}a_{k}\frac{t^{k}}{k!}\right|a_{k}\in\mathbb{C}_{p}\right\} $
be the algebra of formal power series in a single variable $t$. Let
$\mathbb{P}$ be the algebra of polynomials in a single vairable $x$
over $\mathbb{C}_{p}$. We denote the action of the linear functional
$L\in\mathbb{P}^{*}$ on a polynomial $p\left(x\right)$
by $\acl L{p\left(x\right)},$ which is linearly extended as $\acl{cL+c^{\prime}L^{\prime}}{p\left(x\right)}=c\acl L{p\left(x\right)}+c^{\prime}\acl{L^{\prime}}{p\left(x\right)}$,
where $c,c^{\prime}\in\mathbb{C}_{p}$. We define a linear functional
on $\mathbb{P}$ by setting
\begin{equation}
\acl{f\left(t\right)}{x^{n}}=a_{n},\quad\text{for all }n\ge0\text{ and }f\left(t\right)\in\mathcal{F}.\label{eq:7}
\end{equation}

By (\ref{eq:7}), we easily get
\begin{equation}
\acl{t^{k}}{x^{n}}=n!\delta_{n,k},\quad\left(n,k\ge0\right),\label{eq:8}
\end{equation}
where $\delta_{n,k}$ is the Kronecker's symbol (see \cite{key-15}).

For $f_{L}\left(t\right)=\sum_{k=0}^{\infty}\acl L{x^{k}}\frac{t^{k}}{k!}$,
we have $\acl{f_{L}\left(t\right)}{x^{n}}=\acl L{x^{n}}$. The map
$L\mapsto f_{L}\left(t\right)$ is vector space isomorphism from $\mathbb{P}^{*}$
onto $\mathcal{F}$. Henceforth $\mathcal{F}$ denotes both the algebra
of formal power series in $t$ and the vector space of all linear
functionals on $\mathbb{P}$, and so an element $f\left(t\right)$
of $\mathcal{F}$ is thought of as both a formal power series
and a linear functional. We call $\mathcal{F}$ the umbral algebra.
The umbral calculus is the study of umbral algebra.

The order $o\left(f\left(t\right)\right)$ of the non-zero power series
$f\left(t\right)$ is the smallest integer $k$ for which the coefficient
of $t^{k}$ does not vanish (see \cite{key-8,key-15}). If $o\left(f\left(t\right)\right)=1$
(respectively, $o\left(f\left(t\right)\right)=0$), then $f\left(t\right)$
is called a delta (respectively, an invertible) series.

For $o\left(f\left(t\right)\right)=1$ and $o\left(g\left(t\right)\right)=0$,
there exists a unique sequence $s_{n}\left(x\right)$ of polynomials
such that $\acl{g\left(t\right)f\left(t\right)^{k}}{s_{n}\left(x\right)}=n!\delta_{n,k}$$\left(n,k\ge0\right)$.
The sequence $s_{n}\left(x\right)$ is called the Sheffer sequence
for $\left(g\left(t\right),f\left(t\right)\right)$, and we write
$s_{n}\left(x\right)\sim\left(g\left(t\right),f\left(t\right)\right)$
(see \cite{key-15}).

For $f\left(t\right)\in\mathcal{F}$ and $p\left(x\right)\in\mathbb{P}$,
by (\ref{eq:8}), we get
\begin{equation}
\acl{e^{yt}}{p\left(x\right)}=p\left(y\right),\quad\acl{f\left(t\right)g\left(t\right)}{p\left(x\right)}=\acl{g\left(t\right)}{f\left(t\right)p\left(x\right)}=\acl{f\left(t\right)}{g\left(t\right)p\left(x\right)}\label{eq:9}
\end{equation}
and
\begin{equation}
f\left(t\right)=\sum_{k=0}^{\infty}\acl{f\left(t\right)}{x^{k}}\frac{t^{k}}{k!},\quad p\left(x\right)=\sum_{k=0}^{\infty}\acl{t^{k}}{p\left(x\right)}\frac{x^{k}}{k!},\quad\left(\text{see \cite{key-15}}\right).\label{eq:10}
\end{equation}

From (\ref{eq:10}), we note that
\begin{equation}
p^{\left(k\right)}\left(0\right)=\acl{t^{k}}{p\left(x\right)}=\acl 1{p^{\left(k\right)}\left(x\right)},\quad\left(k\ge0\right),\label{eq:11}
\end{equation}
where $p^{\left(k\right)}\left(0\right)$ denotes the $k$-th derivative
of $p\left(x\right)$ with respect to $x$ at $x=0$.

By (\ref{eq:11}), we get
\begin{equation}
t^{k}p\left(x\right)=p^{\left(k\right)}\left(x\right)=\frac{d^{k}}{dx^{k}}p\left(x\right),\quad\left(k\ge0\right).\label{eq:12}
\end{equation}

In \cite{key-15}, it is known that
\begin{equation}
s_{n}\left(x\right)\sim\left(g\left(t\right),f\left(t\right)\right)\iff\frac{1}{g\left(\overline{f}\left(t\right)\right)}e^{x\overline{f}\left(t\right)}=\sum_{n=0}^{\infty}s_{n}\left(x\right)\frac{t^{n}}{n!},\quad\left(x\in\mathbb{C}_{p}\right),\label{eq:13}
\end{equation}
where $\overline{f}\left(t\right)$ is the compositional inverse of
$f\left(t\right)$ such that $f\left(\overline{f}\left(t\right)\right)=\overline{f}\left(f\left(t\right)\right)=t$.

From (\ref{eq:12}), we can easily derive the following equation:
\begin{equation}
e^{yt}p\left(x\right)=p\left(x+y\right),\quad\text{where }p\left(x\right)\in\mathbb{P}=\mathbb{C}_{p}\left[x\right].\label{eq:14}
\end{equation}

In this paper, we study degenerate poly-Bernoulli numbers and polynomials
associated with polylogarithm function and $p$-adic invariant integral
on $\Zp$. Finally, we give some identities of those numbers and polynomials
which are derived from umbral calculus.

\section{Some applications of degenerate poly-Bernoulli numbers}

Now, we consider the degenerate poly-Bernoulli polynomials which are
given by the generating function
\begin{equation}
\frac{\li_{k}\left(1-\left(1+\lambda t\right)^{-\frac{1}{\lambda}}\right)}{\left(1+\lambda t\right)^{\frac{1}{\lambda}}-1}e^{xt}=\sum_{n=0}^{\infty}\beta_{n,\lambda}^{\left(k\right)}\left(x\right)\frac{t^{n}}{n!},\quad\left(k\in\mathbb{Z}\right).\label{eq:15}
\end{equation}

From (\ref{eq:13}) and (\ref{eq:15}), we have
\begin{equation}
\beta_{n,\lambda}^{\left(k\right)}\left(x\right)\sim\left(\frac{\left(1+\lambda t\right)^{\frac{1}{\lambda}}-1}{\li_{k}\left(1-\left(1+\lambda t\right)^{-\frac{1}{\lambda}}\right)},t\right),\label{eq:16}
\end{equation}
and
\begin{equation}
\beta_{n,\lambda}^{\left(k\right)}\left(x\right)=\sum_{l=0}^{n}\binom{n}{l}\beta_{l,\lambda}^{\left(k\right)}x^{n-l},\label{eq:17}
\end{equation}
where $\beta_{l,\lambda}^{\left(k\right)}=\beta_{l,\lambda}^{\left(k\right)}\left(0\right)$
are called the degenerate poly-Bernoulli numbers.

Thus, by (\ref{eq:17}), we get
\begin{align}
\int_{x}^{x+y}\beta_{n,\lambda}^{\left(k\right)}\left(u\right)du & =\frac{1}{n+1}\left\{ \beta_{n+1,\lambda}^{\left(k\right)}\left(x+y\right)-\beta_{n+1,\lambda}^{\left(k\right)}\left(x\right)\right\} \label{eq:18}\\
 & =\frac{e^{yt}-1}{t}\beta_{n,\lambda}^{\left(k\right)}\left(x\right).\nonumber
\end{align}

Let $f\left(t\right)$ be the linear functional such that
\[
\acl{f\left(t\right)}{p\left(x\right)}=\int_{\Zp}\frac{\left(e^{t}-1\right)\li_{k}\left(1-\left(1+\lambda t\right)^{-\frac{1}{\lambda}}\right)}{t\left(\left(1+\lambda t\right)^{\frac{1}{\lambda}}-1\right)}p\left(x\right)d\mu_{0}\left(x\right)
\]
for all polynomials $p\left(x\right)$. Then it can be determined as follows: for any $p(x)\in \mathbb{P}$,
\[ \acl{\frac{t}{e^t -1}}{p(x)} =\int_{\Zp} p(x) d\mu_0 (x).\]
Replacing $p(x)$ by $\frac{e^t -1}{t} h(t)p(x)$, for $h(t)\in \mathcal{F}$, we get
\begin{equation}\label{eq:19}
\acl{h(t)}{p(x)}=\int_{\Zp} \frac{e^t-1}{t} h(t) p(x)d\mu_0 (x).
\end{equation}
In particular, for $h(t)=1$, we obtain
\begin{equation}\label{eq:19-1}
\int_{\Zp} \frac{e^t-1}{t} p(x)d\mu_0 (x)=p(0).
\end{equation}

Therefore, by (\ref{eq:19}) and (\ref{eq:19-1}), we obtain the following theorem as a special case.
\begin{thm}
\label{thm:1} For $p\left(x\right)\in\mathbb{P}$, we have
\begin{align*}
&\relphantom{=}
\acl{\frac{\li_{k}\left(1-\left(1+\lambda t\right)^{-\frac{1}{\lambda}}\right)}{\left(1+\lambda t\right)^{\frac{1}{\lambda}}-1}}{p\left(x\right)}\\
&=\int_{\Zp}\frac{\left(e^{t}-1\right)\li_{k}\left(1-\left(1+\lambda t\right)^{-\frac{1}{\lambda}}\right)}{t\left(\left(1+\lambda t\right)^{\frac{1}{\lambda}}-1\right)}p\left(x\right)d\mu_{0}\left(x\right),
\end{align*}
and
\begin{align*}
&\relphantom{=}\acl{\frac{\left(e^{t}-1\right)\li_{k}\left(1-\left(1+\lambda t\right)^{-\frac{1}{\lambda}}\right)}{t\left(\left(1+\lambda t\right)^{\frac{1}{\lambda}}-1\right)}\int_{\Zp}e^{yt}d\mu_{0}\left(y\right)}{p\left(x\right)}\\
&=\int_{\Zp}\frac{\left(e^{t}-1\right)\li_{k}\left(1-\left(1+\lambda t\right)^{-\frac{1}{\lambda}}\right)}{t\left(\left(1+\lambda t\right)^{\frac{1}{\lambda}}-1\right)}p\left(x\right)d\mu_{0}\left(x\right).
\end{align*}

In particular,
\[
\beta_{n,\lambda}^{\left(k\right)}=\acl{\frac{\left(e^{t}-1\right)\li_{k}\left(1-\left(1+\lambda t\right)^{-\frac{1}{\lambda}}\right)}{t\left(\left(1+\lambda t\right)^{\frac{1}{\lambda}}-1\right)}\int_{\Zp}e^{yt}d\mu_{0}\left(y\right)}{x^{n}},\quad\left(n\ge0\right).
\]
\end{thm}
\bigskip

Note that
\begin{align*}
 & \acl{\int_{\Zp}e^{yt}d\mu_{0}\left(y\right)}{\frac{e^{t}-1}{t}\beta_{n,\lambda}^{\left(k\right)}\left(x\right)}\\
= & \frac{1}{n+1}\acl{\frac{t}{e^{t}-1}}{\beta_{n+1,\lambda}^{\left(k\right)}\left(x+1\right)-\beta_{n+1,\lambda}^{\left(k\right)}\left(x\right)}\\
= & \frac{1}{n+1}\sum_{l=0}^{n+1}\binom{n+1}{l}B_{l}\left(\beta_{n+1-l,\lambda}^{\left(k\right)}\left(1\right)-\beta_{n+1-l,\lambda}^{\left(k\right)}\right)=\beta_{n,\lambda}^{\left(k\right)}.
\end{align*}

It is easy to show that
\begin{align}
 & \frac{\left(e^{t}-1\right)\li_{k}\left(1-\left(1+\lambda t\right)^{-\frac{1}{\lambda}}\right)}{t\left(\left(1+\lambda t\right)^{\frac{1}{\lambda}}-1\right)}\sum_{n=0}^{\infty}\int_{\Zp}\left(x+y\right)^{n}d\mu_{0}\left(y\right)\frac{t^{n}}{n!}\label{eq:20}\\
 & =\frac{\left(e^{t}-1\right)\li_{k}\left(1-\left(1+\lambda t\right)^{-\frac{1}{\lambda}}\right)}{t\left(\left(1+\lambda t\right)^{\frac{1}{\lambda}}-1\right)}\times\frac{t}{e^{t}-1}e^{xt}\nonumber \\
 & =\sum_{n=0}^{\infty}\beta_{n,\lambda}^{\left(k\right)}\left(x\right)\frac{t^{n}}{n!}.\nonumber
\end{align}

Thus, by (\ref{eq:20}), we get
\begin{align}
\beta_{n,\lambda}^{\left(k\right)}\left(x\right) & =\frac{\left(e^{t}-1\right)\li_{k}\left(1-\left(1+\lambda t\right)^{-\frac{1}{\lambda}}\right)}{t\left(\left(1+\lambda t\right)^{\frac{1}{\lambda}}-1\right)}\int_{\Zp}\left(x+y\right)^{n}d\mu_{0}\left(y\right)\label{eq:21}\\
 & =\frac{\li_{k}\left(1-\left(1+\lambda t\right)^{-\frac{1}{\lambda}}\right)}{\left(1+\lambda t\right)^{\frac{1}{\lambda}}-1}x^{n}\nonumber
\end{align}

Therefore, by (\ref{eq:21}), we obtain the following theorem.
\begin{thm}
\label{thm:2} For $p\left(x\right)\in\mathbb{P}$, we have
\begin{align*}
&\relphantom{=}\frac{\left(e^{t}-1\right)\li_{k}\left(1-\left(1+\lambda t\right)^{-\frac{1}{\lambda}}\right)}{t\left(\left(1+\lambda t\right)^{\frac{1}{\lambda}}-1\right)}\int_{\Zp}p\left(x+y\right)d\mu_{0}\left(y\right) \\
& =\frac{\left(e^{t}-1\right)\li_{k}\left(1-\left(1+\lambda t\right)^{-\frac{1}{\lambda}}\right)}{t\left(\left(1+\lambda t\right)^{\frac{1}{\lambda}}-1\right)}\int_{\Zp}e^{yt}p\left(x\right)d\mu_{0}\left(y\right)\\
 & =\frac{\li_{k}\left(1-\left(1+\lambda t\right)^{-\frac{1}{\lambda}}\right)}{\left(1+\lambda t\right)^{\frac{1}{\lambda}}-1}p\left(x\right).
\end{align*}
\end{thm}
\bigskip

For $r\in\mathbb{N}$, let us consider the higher-order degenerate
poly-Bernoulli polynomials as follows:

\begin{align}
 & \left(\frac{\left(e^{t}-1\right)\li_{k}\left(1-\left(1+\lambda t\right)^{-\frac{1}{\lambda}}\right)}{t\left(\left(1+\lambda t\right)^{\frac{1}{\lambda}}-1\right)}\right)^{r}\int_{\Zp}\cdots\int_{\Zp}e^{\left(x_{1}+\cdots+x_{r}+x\right)t}d\mu_{0}\left(x_{1}\right)\cdots d\mu_{0}\left(x_{r}\right)\label{eq:22}\\
= & \left(\frac{\li_{k}\left(1-\left(1+\lambda t\right)^{-\frac{1}{\lambda}}\right)}{\left(1+\lambda t\right)^{\frac{1}{\lambda}}-1}\right)^{r}e^{xt}=\sum_{n=0}^{\infty}\beta_{n,\lambda}^{\left(k,r\right)}\left(x\right)\frac{t^{n}}{n!}.\nonumber
\end{align}
Thus, we obtain
\begin{align}
\beta_{n,\lambda}^{\left(k,r\right)}\left(x\right) & =\left(\frac{\li_{k}\left(1-\left(1+\lambda t\right)^{-\frac{1}{\lambda}}\right)}{\left(1+\lambda t\right)^{\frac{1}{\lambda}}-1}\right)^{r}x^{n}\label{eq:27}\\
 & =\left(\frac{\left(e^{t}-1\right)\li_{k}\left(1-\left(1+\lambda t\right)^{-\frac{1}{\lambda}}\right)}{t\left(\left(1+\lambda t\right)^{\frac{1}{\lambda}}-1\right)}\right)^{r}\nonumber\\
&\relphantom{=}\times\int_{\Zp}\cdots\int_{\Zp}\left(x_{1}+\cdots+x_{r}+x\right)^{n}d\mu_{0}\left(x_{1}\right)\cdots d\mu_{0}\left(x_{r}\right),\nonumber
\end{align}
where $n\ge0$.

Here, for $x=0$, $\beta_{n,\lambda}^{\left(k,r\right)}=\beta_{n,\lambda}^{\left(k,r\right)}\left(0\right)$
are called the  degenerate poly-Bernoulli numbers of order $r$. From
(\ref{eq:22}), we note that
\begin{equation}
\beta_{n,\lambda}^{\left(k\right)}\left(x\right)\sim\left(\left(\frac{\left(1+\lambda t\right)^{\frac{1}{\lambda}}-1}{\li_{k}\left(1-\left(1+\lambda t\right)^{-\frac{1}{\lambda}}\right)}\right)^{r},t\right).\label{eq:23}
\end{equation}

Therefore, by (\ref{eq:27}), we obtain the following
theorem.
\begin{thm}
\label{thm:3} For $p\left(x\right)\in\mathbb{P}$ and $r\in\mathbb{N}$,
we have
\begin{align*}
 &\relphantom{=} \left(\frac{\left(e^{t}-1\right)\li_{k}\left(1-\left(1+\lambda t\right)^{-\frac{1}{\lambda}}\right)}{t\left(\left(1+\lambda t\right)^{\frac{1}{\lambda}}-1\right)}\right)^{r}\int_{\Zp}\cdots\int_{\Zp}p\left(x_{1}+\cdots+x_{r}+x\right)d\mu_{0}\left(x_{1}\right)\cdots d\mu_{0}\left(x_{r}\right)\\
  &= \left(\frac{\left(e^{t}-1\right)\li_{k}\left(1-\left(1+\lambda t\right)^{-\frac{1}{\lambda}}\right)}{t\left(\left(1+\lambda t\right)^{\frac{1}{\lambda}}-1\right)}\right)^{r}\int_{\Zp}\cdots\int_{\Zp}e^{(x_1+\cdots+x_r)t}p(x)d\mu_{0}\left(x_{1}\right)\cdots d\mu_{0}\left(x_{r}\right)\\
 &= \left(\frac{\li_{k}\left(1-\left(1+\lambda t\right)^{\frac{1}{\lambda}}\right)}{\left(1+\lambda t\right)^{\frac{1}{\lambda}}-1}\right)^{r}p\left(x\right).
\end{align*}
\end{thm}
\bigskip
Let us consider the linear functional $f_{r}\left(t\right)$ such that
\begin{align}\label{eq:28}
&\relphantom{=}\acl{f_{r}\left(t\right)}{p\left(x\right)}\\
&=\int_{\Zp}\cdots\int_{\Zp}\left(\frac{\left(e^{t}-1\right)\li_{k}\left(1-\left(1+\lambda t\right)^{-\frac{1}{\lambda}}\right)}{t\left(\left(1+\lambda t\right)^{\frac{1}{\lambda}}-1\right)}\right)^{r}\left.p(x)\right|_{x=x_1+\cdots+x_r}d\mu_{0}\left(x_{1}\right)\cdots d\mu_{0}\left(x_{r}\right)\nonumber
\end{align}
for all polynomials $p\left(x\right)$. Then it can be determined in the following way: for $p(x)\in\mathbb{P}$,
\begin{align*}
&\acl{\left(\frac{t}{e^t-1}\right)^r}{p(x)}=\int_{\Zp}\cdots\int_{\Zp} \left.p(x)\right|_{x=x_1+\cdots+x_r} d\mu_0(x_1)\cdots d\mu_0 (x_r).
\end{align*}
Replacing $p(x)$ by $\left(\frac{e^t-1}{t}h(t)\right)^rp(x)$, for $h(t)\in\mathcal{F}$, we have
\begin{equation}\label{Star-1}
\acl{h(t)^r}{p(x)}=\int_{\Zp}\cdots\int_{\Zp} \left(\frac{e^t-1}{t} h(t)\right)^r \left.p(x)\right|_{x=x_1+\cdots+x_r} d\mu_0 (x_1)\cdots d\mu_0 (x_r).
\end{equation}
In particular, for $h(t)=1$, we get
\begin{equation}\label{Star-2}
\int_{\Zp}\cdots\int_{\Zp} \left(\frac{e^t-1}{t}\right)^r \left.p(x)\right|_{x=x_1+\cdots+x_r}d\mu_0 (x_1)\cdots d\mu_0 (x_r) =p(0).
\end{equation}

Therefore, by (\ref{Star-1}) and (\ref{Star-2}), we obtain the following
theorem.
\begin{thm}
\label{thm:4} For $p\left(x\right)\in\mathbb{P}$, we have
\begin{align*}
 & \acl{\left(\frac{\li_{k}\left(1-\left(1+\lambda t\right)^{-\frac{1}{\lambda}}\right)}{\left(1+\lambda t\right)^{\frac{1}{\lambda}}-1}\right)^{r}}{p\left(x\right)}\\
= & \int_{\Zp}\cdots\int_{\Zp}\left(\frac{\left(e^{t}-1\right)\li_{k}\left(1-\left(1+\lambda t\right)^{-\frac{1}{\lambda}}\right)}{t\left(\left(1+\lambda t\right)^{\frac{1}{\lambda}}-1\right)}\right)^{r}\left.p(x)\right|_{x=x_1+\cdots+x_r}d\mu_{0}\left(x_{1}\right)\cdots d\mu_{0}\left(x_{r}\right),
\end{align*}
and
\begin{align*}
 & \acl{\left(\frac{\left(e^{t}-1\right)\li_{k}\left(1-\left(1+\lambda t\right)^{-\frac{1}{\lambda}}\right)}{t\left(\left(1+\lambda t\right)^{\frac{1}{\lambda}}-1\right)}\right)^{r}\int_{\Zp}\cdots\int_{\Zp}e^{\left(x_{1}+\cdots+x_{r}\right)t}d\mu_{0}\left(x_{1}\right)\cdots d\mu_{0}\left(x_{r}\right)}{p\left(x\right)}\\
= &\int_{\Zp}\cdots\int_{\Zp} \left(\frac{\left(e^{t}-1\right)\li_{k}\left(1-\left(1+\lambda t\right)^{-\frac{1}{\lambda}}\right)}{t\left(\left(1+\lambda t\right)^{\frac{1}{\lambda}}-1\right)}\right)^{r}\left.p(x)\right|_{x=x_1+\cdots+x_r}d\mu_{0}\left(x_{1}\right)\cdots d\mu_{0}\left(x_{r}\right).
\end{align*}

In particular,
\begin{align*}
\beta_{n,\lambda}^{\left(k,r\right)} & =\acl{\left(\frac{\left(e^{t}-1\right)\li_{k}\left(1-\left(1+\lambda t\right)^{-\frac{1}{\lambda}}\right)}{t\left(\left(1+\lambda t\right)^{\frac{1}{\lambda}}-1\right)}\right)^{r}\int_{\Zp}\cdots\int_{\Zp}e^{\left(x_{1}+\cdots+x_{r}\right)t}d\mu_{0}\left(x_{1}\right)\cdots d\mu_{0}\left(x_{r}\right)}{x^{n}}.
\end{align*}
\end{thm}
\bigskip

\begin{rmk}
It is not difficult to show that
\begin{align*}
 & \acl{\left(\frac{\left(e^{t}-1\right)\li_{k}\left(1-\left(1+\lambda t\right)^{-\frac{1}{\lambda}}\right)}{t\left(\left(1+\lambda t\right)^{\frac{1}{\lambda}}-1\right)}\right)^{r}\int_{\Zp}\cdots\int_{\Zp}e^{\left(x_{1}+\cdots+x_{r}\right)t}d\mu_{0}\left(x_{1}\right)\cdots d\mu_{0}\left(x_{r}\right)}{x^{n}}\\
= & \sum_{n=n_{1}+\cdots+n_{r}}\binom{n}{n_{1},\dots,n_{r}}\acl{\frac{\left(e^{t}-1\right)\li_{k}\left(1-\left(1+\lambda t\right)^{-\frac{1}{\lambda}}\right)}{t\left(\left(1+\lambda t\right)^{\frac{1}{\lambda}}-1\right)}\int_{\Zp}e^{x_{n_{1}}t}d\mu_{0}\left(x_{1}\right)}{x^{m_{1}}}\times\cdots\\
 & \times\acl{\frac{\left(e^{t}-1\right)\li_{k}\left(1-\left(1+\lambda t\right)^{-\frac{1}{\lambda}}\right)}{t\left(\left(1+\lambda t\right)^{\frac{1}{\lambda}}-1\right)}\int_{\Zp}e^{x_{n_{r}}t}d\mu_{0}\left(x_{n_{r}}\right)}{x^{n_{r}}}.
\end{align*}

Thus, we get
\[
\beta_{n,\lambda}^{\left(k,r\right)}=\sum_{n=n_{1}+\cdots+n_{r}}\binom{n}{n_{1},\dots,n_{r}}\beta_{n_{1},\lambda}^{\left(k\right)}\cdots\beta_{n_{r},\lambda}^{\left(k\right)}.
\]
\end{rmk}

\bibliographystyle{amsplain}
\providecommand{\bysame}{\leavevmode\hbox to3em{\hrulefill}\thinspace}
\providecommand{\MR}{\relax\ifhmode\unskip\space\fi MR }
\providecommand{\MRhref}[2]{%
  \href{http://www.ams.org/mathscinet-getitem?mr=#1}{#2}
}
\providecommand{\href}[2]{#2}

\end{document}